\documentclass[11pt]{article}
\usepackage{amsmath}
\usepackage{amsthm}
\usepackage{amssymb}

\newtheorem{theo}{Theorem}[section]

\newcommand{\FF}{{\cal F}}
\newcommand{\GG}{{\cal G}}

\begin{document}
\date{}

\title{
Partitioning all $k$-subsets into $r$-wise intersecting families
}

\author{Noga Alon
\thanks
{Princeton University,
Princeton, NJ 08544, USA and
Tel Aviv University, Tel Aviv 69978,
Israel.
Email: {\tt nalon@math.princeton.edu}.  
Research supported in part by
NSF grant DMS-1855464 and BSF grant 2018267.}
}

\maketitle
\begin{abstract}
Let $r \geq 2$, $n$ and $k$ be integers satisfying
$k \leq \frac{r-1}{r}n$. In the original arXiv version of this
note we suggested a  conjecture that
the family of all $k$-subsets of an $n$-set
cannot be partitioned into fewer
than  $\lceil n-\frac{r}{r-1}(k-1) \rceil$ $r$-wise intersecting
families. We noted that if true this is tight for all values of the parameters,
that the case $r=2$ is Kneser's conjecture, proved by
Lov\'asz, and observed that the assertion also holds provided
$r$ is either a prime number or a power of $2$. We have recently learned,
however, that the assertion of the conjecture for all values of the
parameters follows from a recent result of Azarpendar and Jafari
\cite{AJ}. 
\end{abstract}

\section{Introduction}
One of the earliest and possibly the best known 
application of topological methods in
extremal combinatorics is the Kneser conjecture, now 
Lov\'asz' Theorem \cite{Lo}. This theorem  asserts that for every
$n \geq 2k$ it is impossible to split the family of all $k$-subsets
of an $n$-set into fewer than $n-2k+2$ intersecting families. 
The main purpose of this brief note is to study
the following extension of this result. Call a family of
subsets $r$-wise intersecting if any collection of at most $r$
subsets in it has a common point.  
\begin{theo}
\label{t11}
Let $r \geq 2, n$ and $k$ be integers, and suppose
$k \leq \frac{r-1}{r}n$. Then the family of all $k$-subsets of 
$[n]=\{1,2, \ldots ,n\}$ cannot be partitioned into fewer
than  $\lceil n-\frac{r}{r-1}(k-1) \rceil$ $r$-wise intersecting
families $\FF_i$. This is tight for all admissible values of the
parameters.
\end{theo}
\vspace{0.1cm}

\noindent
The case $r=2$ is Kneser's Conjecture proved by Lov\'asz.
\vspace{0.1cm}

\noindent
As mentioned in the abstract, we have learned that
this theorem is in fact a special
case of a recent result of 
Azarpendar and Jafari, see Theorem 1.1 and Corollary 1.1 
in \cite{AJ}. Since we have already posted an early version of
this note before learning about this paper, we still describe
here the observation that
the assertion of the theorem holds 
for every prime $r$ and for every $r$ which is a power of $2$. 

\section{Proof}
One can prove Theorem \ref{t11} for prime $r$ and for $r$ which is
a power of $2$ using tools from equivariant
topology. It is much simpler, however, to deduce it from known
results about the chromatic numbers of stable Kneser hypergraphs.
We start by describing these results.

For $n \geq rk$ let
$KG^r(k,n)$ denote the Kneser hypergraph whose vertex set is the
set of all ${n \choose k}$ subsets of cardinality $k$ of an $n$-set
$[n]=\{1,2,\ldots ,n\}$ and whose edges are all $r$-tuples of
$k$-sets which are pairwise disjoint. For $r=2$ this is the Kneser
graph. Lov\'asz proved that the chromatic number of $KG^{2}(k,n)$
is $n-2k+2$ for all $n \geq 2k$, and in \cite{AFL} it is proved that
the chromatic number of $KG^r(k,n)$ is 
$\lceil \frac{n-r(k-1)}{r-1} \rceil$ for all $n \geq rk$.

Call a subset $F$ of $[n]$ $s$-stable if any two elements of $F$
are at distance at least $s$ in the
cyclic order on $[n]$. Let $KG^{r}(k,n)_{s-stab}$ be the induced
sub-hypergraph of $KG^{r}(k,n) $ on the set of vertices which are
$s$-stable. Ziegler \cite{Zi},
and Drewnowski, {\L}uczak and the present author \cite{ADL}
conjectured  that the chromatic number of $KG^{r}(k,n)_{r-stab}$
is also $\lceil \frac{n-r(k-1)}{r-1} \rceil$, just as that
of the full hypergraph $KG^{r}(k,n)$.  
We need the following known result regarding this conjecture.
\begin{theo}[\cite{Sch}, \cite{ADL}]
\label{t21}
If $r$ is any power of $2$ then the
chromatic number of $KG^{r}(k,n)_{r-stab}$ 
is $\lceil \frac{n-r(k-1)}{r-1} \rceil$ 
for all admissible values of $k,n$.
\end{theo} 
The case $r=2$ was proved by Schrijver in \cite{Sch},  and in
\cite{ADL} it is shown that if the result holds for $r_1$ and $r_2$
then it also holds for $r_1r_2$, implying the assertion for all
powers of $2$.

Improving results of  Meunier \cite{Me} and of
Alishahi and Hajiabolhassan \cite{AH}, Frick proved
in \cite{Fr} that the chromatic number of
$KG^{r}(k,n)_{2-stab}$ is $\lceil \frac{n-r(k-1)}{r-1} \rceil$. 
Another result proved in \cite{Fr} is (a slightly stronger
version of)
the following.
\begin{theo}[\cite{Fr}, Theorem 3.10]
\label{t22}
For any prime $r$ and for any 
partition of $[n]$
into subsets $C_i$, each of size at most $r-1$, the induced
subhypergraph of $KG^{r}(k,n)$ on the set of all vertices $F$
that contain at most $1$ element of each $C_i$  has chromatic
number $\lceil \frac{n-r(k-1)}{r-1} \rceil$.
\end{theo}
\vspace{0.2cm}

\noindent
{\bf Proof of Theorem \ref{t11} for $r$ prime or a power of $2$:}\, 
The upper bound showing that the result is tight 
(for all values of $r$) is simple.
Put $s=\lfloor \frac{rk-1}{r-1} \rfloor$. Then every collection
of $r$ subsets of cardinality $k$ of the subset
$S=\{n-s+1,n-s+2, \ldots ,n\}$ has a common point. For every
$i \leq n-s$ let $F_{i}$ be the collection of all $k$-subsets of $[n]$
whose minimum element is $i$. Let $F_{n-s+1}$ be the remaining
$k$-subsets, that is, all those contained in $S$. This shows
that there is a construction with
$$
m =n-s+1 =n -\lfloor \frac{rk-1}{r-1} \rfloor+1=
\lceil n-\frac{r}{r-1}(k-1) \rceil,
$$
as needed.

To prove the lower bound let $n,k,r$ and $\FF_1, \ldots \FF_m$
be as in the  statement of the theorem. Assume, first, that 
$r$ is a prime. Let $C_1, C_2, \ldots ,C_n$ be $n$ pairwise
disjoint sets, each of size $r-1$, and let $C=\cup_{i=1}^n C_i$
be their union. For each of the families 
$\FF_i$, let $\GG_i$ be the family of all $k$-subsets of $C$
obtained as follows. For each member $F=\{i_1, i_2, \ldots ,i_k\}$
of $\FF_i$, let $C(F)$ denote the family of all
$(r-1)^k$ subsets of $C$ containing exactly one element
of each $C_{i_j}$ for $1 \leq j \leq k$. The family $\GG_i$ is the
union of all families $C(F)$ for $F \in \FF_i$.

We claim that
no set $\GG_i$ contains $r$ pairwise disjoint sets. Indeed,
every collection of $r$ (not necessarily  distinct)
members $G_1,G_2, \ldots ,G_r$ of $\GG_i$ consists of
subsets that belong to $C(F_1), C(F_2), \ldots , C(F_r)$,
respectively, for some
(not necessarily  distinct) members $F_j \in \FF_i$. Since
$\FF_i$ is $r$-wise intersecting there is a common point, say
$\ell$,
in all sets  $F_j$. Thus each $G_j$ contains a point of
$C_{\ell}$ and as $|C_{\ell}|=r-1$ some pair of sets $G_j$
intersect inside $C_{\ell}$, by the pigeonhole principle.
This proves the claim.  

Note that the union of all families
$\GG_i$ is exactly the collection of 
all the $k$-subsets of $C$ that contain at most
$1$ element from each $C_i$. Therefore the families $\GG_i$ provide
a proper coloring of the hypergraph described in Theorem
\ref{t22}  with parameters $(r-1)n$, $r$ and $k$. The chromatic
number of this hypergraph is 
$\lceil \frac{(r-1)n-r(k-1)}{r-1} \rceil$, providing the required
lower bound for $m$.

The proof for $r$ which is a power of $2$ is similar, using
the result in Theorem \ref{t21}. We apply the same construction
with sets $C_i$ and families $\GG_i$ as before, and place the sets
$C_i$ along a cycle of length $(r-1)n$, where each set $C_i$
appears contiguously along the cycle. It is then easy to see that
the union of the families $\GG_i$ contains a copy of the
hypergraph $KG^{r}((r-1)n,k)_{r-stab}$ on the set of vertices
$C=\cup_{i=1}^n C_i$ of the cycle (as well as some additional hyperedges). 
As before, here too each family $\GG_i$ contains
no $r$ pairwise disjoint members. 
The lower bound thus follows
from Theorem \ref{t21}.  This completes the proof.
\hfill $\Box$
\vspace{0.2cm}

\noindent
{\bf Acknowledgment}
I thank Ryan Alweiss for many helpful discussions and Amir Jafari
for telling me about \cite{AJ}. 

\end{document}